\documentclass[10pt]{article}

\newcommand{\update}{2024-11-11} 

\title{\vskip-1.0em\bf\newheadfont{Small values and forbidden values for the Fourier anti-diagonal constant of a finite group}}
\author{\sc Yemon Choi}
\date{\sc \update}

\usepackage{amsmath,amssymb,amsfonts}

\newcommand{\newheadfont}{\bfseries} 

\usepackage[a4paper,left=29mm,right=29mm,top=29mm,bottom=29mm,footskip=12mm]{geometry}

\usepackage{mathpazo} 
\usepackage{sectsty}
\allsectionsfont{\scshape} 

\renewcommand{\newheadfont}{\bfseries\scshape}

\usepackage[colorlinks=true,linkcolor=blue,citecolor=red]{hyperref}
\usepackage[dvipsnames]{xcolor}

\newcommand{\dt}[1]{{\bfseries\itshape\textcolor{Bittersweet}{#1}\/}} 

\usepackage{amsthm}
\newcounter{pulse}[section]
\numberwithin{pulse}{section}

\newcommand{\thmsep}{\topsep} 

\newtheoremstyle{newplain} 
  {\thmsep}   
  {\thmsep}   
  {\itshape}  
  {0pt}       
  {\newheadfont} 
  {.}         
  {5pt plus 1pt minus 1pt} 
  {}          

\newtheoremstyle{newdef} 
  {\thmsep}   
  {\thmsep}   
  {\normalfont}  
  {0pt}       
  {\newheadfont} 
  {.}         
  {5pt plus 1pt minus 1pt} 
  {}          

\newtheoremstyle{newrem}
  {\thmsep}   
  {\thmsep}   
  {\normalfont}  
  {0pt}       
  {\newheadfont} 
  {.}         
  {5pt plus 1pt minus 1pt} 
  {}          

\theoremstyle{newplain}
\newtheorem{thm}[pulse]{Theorem}
\newtheorem{prop}[pulse]{Proposition}
\newtheorem{lem}[pulse]{Lemma}
\newtheorem{cor}[pulse]{Corollary}
\theoremstyle{newdef}

\newtheorem{eg}[pulse]{Example}
\theoremstyle{newrem}
\newtheorem{rem}[pulse]{Remark}

\numberwithin{equation}{section}
\newcommand{\defeq}{:=}
\renewcommand{\emph}[1]{{\bfseries#1}}
\newcommand{\veps}{\varepsilon}

\newenvironment{romnum}{%
\begin{enumerate}

}{\end{enumerate}\ignorespacesafterend}

\newenvironment{alphnum}{%
\begin{enumerate}

}{\end{enumerate}\ignorespacesafterend}

\newcommand{\Cplx}{\mathbb C}
\newcommand{\Nat}{\mathbb N}

\newcommand{\bbF}{{\mathbb F}}
\newcommand{\bbT}{{\mathbb T}}
\newcommand{\cL}{{\mathcal L}}

\newcommand{\gpid}[1][]{{\mathsf{1}_{#1}}} 
\renewcommand{\gpid}[1][]{e_{#1}} 
\newcommand{\cd}[1]{{#1}(\gpid)} 
\renewcommand{\cd}[1]{d({#1})} 
\newcommand{\abs}[1]{{\lvert #1 \rvert}} 
\newcommand{\gpind}[2]{\abs{#1:#2}} 
\DeclareMathOperator{\AD}{AD}
\DeclareMathOperator{\Irr}{Irr}
\DeclareMathOperator{\maxdeg}{maxdeg} 
\DeclareMathOperator{\mindeg}{mindeg}
\DeclareMathOperator{\cp}{cp}
\DeclareMathOperator{\cdset}{cd}

\newcommand{\Ind}{{\operatorname{Ind}\nolimits}} 
\newcommand{\Res}[2]{{#1\rvert}_{#2}}
\newcommand{\ip}[2]{{\langle {#1}, {#2}\rangle}}
\newcommand{\FA}{\operatorname{A}}  
\DeclareMathOperator{\Stab}{Stab} 
\DeclareMathOperator{\supp}{supp}
\DeclareMathOperator{\Tr}{Tr} 

\newcommand{\SL}{\mathrm{SL}}
\newcommand{\SO}{\mathrm{SO}}
\newcommand{\SU}{\mathrm{SU}}
\newcommand{\PSL}{\mathrm{PSL}}
\newcommand{\GL}{\mathrm{GL}}
\newcommand{\PGL}{\mathrm{PGL}}
\newcommand{\PSp}{\mathrm{PSp}}

\begin{document}
\maketitle

\begin{abstract}
For a finite group $G$, let $\AD(G)$ denote the Fourier norm of the antidiagonal in $G\times G$. It was shown recently in \cite{YC_minorant} that $\AD(G)$ coincides with the amenability constant of the Fourier algebra of $G$, and is equal to the normalized sum of the cubes of character degrees of $G$. Motivated by a gap result for amenability constants from~\cite{BEJ_AG}, we determine exactly which numbers in the interval $[1,2]$ arise as values of $\AD(G)$.

As a by-product, we show that the set of values of $\AD(G)$ does not contain all its limit points. Some other calculations or bounds for $\AD(G)$ are given for familiar classes of finite groups. We also indicate a connection between $\AD(G)$ and the commuting probability of $G$, and use this to show that every finite group $G$ satisfying $\AD(G)< \frac{61}{15}$ must be solvable; here the value $\frac{61}{15}$ is best possible.


\end{abstract}


\begin{section}{Introduction}

\begin{subsection}{Background context}
Given a finite group $G$, the algebra of complex-valued functions on $G$ (equipped with pointwise product) only depends on the cardinality of $G$, and does not detect the group structure. However, there is a canonical submultiplicative norm on this algebra, the \dt{Fourier norm}, such that the resulting normed algebra $\FA(G)$ characterizes the starting group $G$ up to isomorphism. (More precisely: given finite groups $G$ and $H$, there is an isometric algebra isomorphism between $\FA(G)$ and $\FA(H)$ if and only if $G$ and $H$ are isomorphic groups; this is a a special case of a result of Walter~\cite{Walter_JFA72}.)

By identifying a subset of $G$ with its indicator function, one can speak of the Fourier norm of a subset of~$G$. Calculating Fourier norms of arbitrary subsets is hard (see \cite{Sanders_CDT-AG} for a systematic approach), but there is one case where an exact calculation is possible and gives interesting answers. Consider the set $\{ (g,g^{-1})\colon g\in G\}$. The Fourier norm of this subset of $G\times G$, denoted by $\AD(G)$ in this paper, is the \dt{Fourier anti-diagonal constant} mentioned in the title.
It was recently shown by the author \cite[Theorem 1.4]{YC_minorant} that we have the following explicit formula for~$\AD(G)$:
\begin{equation}\label{eq:formula for AD}
\AD(G) = \frac{1}{\abs{G}}\sum_{\varphi\in\Irr(G)} \cd{\varphi}^3
\end{equation}
where $\Irr(G)$ is the set of irreducible complex characters of $G$ and $\cd{\varphi}$ is the degree of~$\varphi$.

Equation \eqref{eq:formula for AD} implies that $\AD(G\times H)=\AD(G)\AD(H)$, and that $\AD$ is invariant under isoclinism. It can also be shown that $\AD(H)\leq\AD(G)$ whenever $H\leq G$ (see Proposition \ref{p:monotone} below). These hereditary properties suggest that $\AD(G)$, viewed as a numerical invariant of $G$, deserves further study.
Furthermore, the sum on the right-hand side of Equation \eqref{eq:formula for AD} had already arisen in earlier work of Johnson \cite{BEJ_AG} on Fourier algebras of compact groups. The results in \cite[\S4]{BEJ_AG} provide an attractive application of the character theory of finite groups to obtain new (counter-)examples in functional analysis. (For a fuller discussion, see \cite[Section 1]{YC_minorant}.)

The following observations, taken from \cite[Proposition 4.3]{BEJ_AG}, are easy consequences of~\eqref{eq:formula for AD}:
\begin{itemize}
\item if $G$ is abelian, then $\AD(G)=1$;
\item if $G$ is non-abelian, then $\AD(G)\geq \frac{3}{2}$.
\end{itemize}
Since $\AD(G^n) = \AD(G)^n$, this shows that $\AD(G)$ can take arbitrarily large values. However, to the author's knowledge, nothing further has been done to study the possible values of $\AD(G)$ as $G$ ranges over all non-abelian finite groups. The purpose of the present paper is to make a start on filling this gap.
\end{subsection}

\begin{subsection}{Our main new results}
The following result has probably been noticed independently by many readers of Johnson's paper, although it is not stated explicitly there. (A proof will be given in Section \ref{s:easy lower bounds} for sake of completeness.)

\begin{prop}[implicitly folklore]\label{p:if maxdeg leq 2}
Let $G$ be a finite group, and suppose that $\cd{\varphi}\leq 2$ for all $\varphi\in\Irr(G)$. Then $\AD(G) \in \{ 2- n^{-1} \colon n\in\Nat\}$.
\end{prop}

Moreover, every number in $\{2-n^{-1} \colon n\in\Nat\}$ is realized as the $\AD$-constant of some (non-unique) finite group: this can be seen by considering cyclic groups and dihedral groups. Our first main result is that these are the \emph{only} values of $\AD$ attained by finite groups in the interval $[1,2]$. To be precise:

\begin{thm}[{The possible values of $\AD(G)$ in $[1,2]$}]
\label{t:values in [1,2]}
Let $G$ be a finite group, and suppose that $\AD(G)\leq 2$. Then
$\AD(G) \in \{2 -n^{-1} \colon n\in\Nat\}$.
\end{thm}

\begin{cor}
The set $\{\AD(G) \colon \text{$G$ a finite group}\}$ is not a closed subset of $[1,\infty)$.
\end{cor}

Theorem \ref{t:values in [1,2]} is an immediate consequence of combining Proposition \ref{p:if maxdeg leq 2} with the following lower bound for $\AD(G)$, which appears to be new.

\begin{prop}\label{p:new lower bound}
Let $G$ be a finite group. If there exists $\varphi\in\Irr(G)$ with $\cd{\varphi}\geq 3$, then $\AD(G)\geq 2+ \abs{G'}^{-1}$.
\end{prop}

The proof of Proposition \ref{p:new lower bound} requires some basic character theory, but nothing harder than Frobenius reciprocity. Perhaps surprisingly, while we do need character theory for finite groups, we do not rely on any structure theory (we do not even need the Sylow theorems).
In contrast, our other main result requires the classification of finite simple groups with characters of small degree.

\begin{thm}[A threshold ensuring solvability]\label{t:implies solvable}
Let $G$ be a finite group. If $\AD(G) < \frac{61}{15}$, then $G$ is solvable.
\end{thm}

A direct calculation shows that $\AD(A_5)= \frac{61}{15}$, and so in this sense Theorem \ref{t:implies solvable} is sharp.
Particular properties of $A_5$, such as its subgroup structure and its Schur multiplier, play an important role in the proof of Theorem \ref{t:implies solvable}, since we need to analyze perfect groups that quotient onto $A_5$.

One difficulty in proving Theorem \ref{t:implies solvable} is that $\AD$ is not monotone (in either direction) with respect to taking quotients, and so knowing that a group $H$ quotients onto $A_5$ does not immediately imply that $\AD(H)\geq\AD(A_5)$.
Instead we require a detour through the \dt{commuting probability} $\cp(G)$ (see the start of Section \ref{s:nonsolvable} for its definition).
Our strategy is inspired by an argument of Tong-Viet in~\cite{TongViet_large}, and indeed the main work needed to prove Theorem~\ref{t:implies solvable} lies in establishing the following stronger version of \cite[Lemma 2.4]{TongViet_large}.

\begin{prop}\label{p:sharper TV}
Let $H$ be a finite non-trivial perfect group satisfying $\cp(H)> \frac{1}{20}$.
Then $H\cong A_5$ or $H\cong\SL(2,5)$.
\end{prop}
\end{subsection}

\begin{subsection}{Outline of this paper}
After some preliminary results in Section~\ref{s:easy lower bounds}, the proof of Proposition~\ref{p:new lower bound} is given in Section~\ref{s:proof of main bound}.
Since the paper is intended for a general audience, we have spelled things out in more detail than specialists in group theory would require.
In Section~\ref{ss:particular} we calculate the values of $\AD(G)$ for some particular families of groups, some of which are related to calculations in earlier sections; and in Section \ref{ss:miscell} we present some partial results on the general theme that ``small values'' of $\AD(G)$ imply that $G$ is close to abelian in some sense. Section \ref{s:nonsolvable} is dedicated to the proof of Proposition \ref{p:sharper TV} and Theorem \ref{t:implies solvable}; this is the only part of the paper that makes use of the theory of the $\cp$ invariant. In the appendix, we collect some proofs of results that are used in the main body of the paper; these results are special cases or weaker versions of known results, but we take the opportunity to provide some extra details and give more elementary arguments.

We finish this introduction by establishing some conventions and fixing notation.
To~reduce unnecessary repetition, we adopt the following convention: henceforth, \emph{all groups are assumed to be finite} unless explicitly stated otherwise.
The identity element of a group $G$ is denoted by~$\gpid$, or $\gpid[G]$ if we wish to avoid ambiguity, and the \dt{derived subgroup} of $G$ (also known as its commutator subgroup) is denoted by~$G'$.

Throughout this article, all representations and characters are taken over the complex field.
The basic representation theory and character theory that we need can be found in several introductory texts, such as~\cite{JamLie}.
We denote the \dt{degree} of a character $\varphi$ by $\cd{\varphi}$; note that this is equal to the dimension of any representation whose trace is $\varphi$.

The set of irreducible characters of $G$ is denoted by $\Irr(G)$, and we write
 $\cdset(G)$ for the set $\{ \cd{\varphi} \colon \varphi\in\Irr(G)\}$ (note that here we are not counting the multiplicities of the irreducible character degrees).
We write $\Irr_n(G)$ for the set of all $\varphi\in\Irr(G)$ that have degree~$n$.
If $G$ is non-abelian, we define $\mindeg(G)\defeq \min\{d\geq 2 \colon \Irr_d(G)\text{ is non-empty}\}$.
For any $G$ (possibly abelian), we define $\maxdeg(G)\defeq \max\{ \cd{\varphi} \colon \varphi\in\Irr(G)\}$.

Finally, given a group $G$, we equip $\Cplx^G$ with the following inner product:
\[
\ip{\varphi}{\psi} \defeq \frac{1}{\abs{G}} \sum_{x\in G} \varphi(x)\overline{\psi(x)} \;.
\]
If $\psi$ is a character of~$G$, then it is irreducible if and only if $\ip{\psi}{\psi}=1$~\cite[Theorem 14.20]{JamLie}.

\end{subsection}

\end{section}

\begin{section}{Some easy lower bounds on $\AD$}\label{s:easy lower bounds}

We start by giving a proof of Proposition \ref{p:if maxdeg leq 2}, since it will also serve as a prototype for later arguments. No novelty is claimed.

\begin{proof}[Proof of Proposition \ref{p:if maxdeg leq 2}]
Since $\cdset(G)\subseteq\{1,2\}$,
\[ \AD(G) = \frac{1}{\abs{G}} \left( \abs{\Irr_1(G)} + 8 \abs{\Irr_2(G)} \right)\,. \]
On the other other hand, basic character theory tells us that
\[ 
  1 = \frac{1}{\abs{G}} \sum_{\varphi\in\Irr(G)} \cd{\varphi}^2 = \frac{1}{\abs{G}} \left( \abs{\Irr_1(G)} + 4 \abs{\Irr_2(G)} \right)\,, \]
and therefore $\AD(G)-2 =  - \abs{\Irr_1(G)}\, \abs{G}^{-1} $.

It is also standard (see e.g.~\cite[Theorem 17.11]{JamLie}) that,
since $\Irr_1(G)$ can be identified with the (Pontrjagin) dual of the abelian group $G/G'$, we have $\abs{\Irr_1(G)} =\gpind{G}{G'}$.
Hence $\AD(G) = 2 - \abs{G'}^{-1}$, and since $\abs{G'}\in\Nat$, the result follows.
\end{proof}

\begin{eg}[Dihedral groups]\label{eg:dihedral}
Let $G$ be a dihedral group of order $2k$, so that $\cdset(G)=\{1,2\}$. If $k$ is odd then $\abs{G'}=k$, and if $k$ is even then $\abs{G'}=\frac{k}{2}$. By repeating the calculation in the proof of Proposition \ref{p:if maxdeg leq 2}, we see that $\AD(G)=2-\frac{1}{k}$ when $k$ is odd and $\AD(G)=2-\frac{2}{k}$ when $k$ is even.
\end{eg}

For general non-abelian $G$, the proof of Proposition \ref{p:if maxdeg leq 2} still suggests a way to proceed. Informally: since $\cd{\varphi}^3 \geq \mindeg(G)\cd{\varphi}^2$ for all $\varphi\in\Irr(G)\setminus \Irr_1(G)$, we can add a correction factor to $\AD(G)$ to obtain something bounded below by $\mindeg(G)$, and the size of the correction factor is controlled by the size of $\abs{G'}$. Making this precise leads us to the following lemma, which provides a convenient tool for dealing with ``generic'' cases.

\begin{lem}[An all-purpose lower bound]\label{l:hammer}
Let $G$ be non-abelian, and let $m,n\in\Nat$ satisfy $\mindeg(G)\geq m$ and $\abs{G'}\geq n$. Then
\begin{equation}
\AD(G) \geq 1+ (m-1)\left(1-\frac{1}{n}\right)\,.
\end{equation} 
\end{lem}

\begin{proof}
Since $\Irr_j(G)$ is empty whenever $2\leq j\leq m-1$, we have
\[
   \AD(G) - \frac{\abs{\Irr_1(G)}}{\abs{G}}
 = \frac{1}{\abs{G}}\sum_{n\geq m} n^3 \abs{\Irr_n(G)} 
\]
and
\[
   1 - \frac{\abs{\Irr_1(G)}}{\abs{G}}
 = \frac{1}{\abs{G}}\sum_{n\geq m} n^2 \abs{\Irr_n(G)}\,.
\]
Hence
\[
\AD(G) - \frac{\abs{\Irr_1(G)}}{\abs{G}} \geq m    \left(1 - \frac{\abs{\Irr_1(G)}}{\abs{G}} \right)\,.
\]
As in the proof of Proposition \ref{p:if maxdeg leq 2},  $\abs{\Irr_1(G)} = \gpind{G}{G'}$. Hence $\abs{\Irr_1(G)} \abs{G}^{-1} \leq n^{-1}$, and plugging this into the previous inequality gives
\[
\AD(G) \geq  m - (m-1) \frac{\abs{\Irr_1(G)}}{\abs{G}} \geq m- \frac{m-1}{n}\,,
\]
completing the proof.
\end{proof}


\begin{cor}[{A sharper form of \cite[Proposition 4.3]{BEJ_AG}}]\label{c:minimal AD}
Let $G$ be non-abelian. Then either $\AD(G) \geq \frac{5}{3}$,
or $\cdset(G)=\{1,2\}$ and $\abs{G'}=2$; in the latter case $\AD(G)=\frac{3}{2}$.
\end{cor}

\begin{proof}
Note that $\mindeg(G)\geq2 \iff \text{$G$ is non-abelian} \iff \abs{G'}\geq 2$.
Therefore, if either $\mindeg(G)\geq 3$ or $\abs{G'}\geq 3$,
applying Lemma \ref{l:hammer} with $(m,n)=(3,2)$ and $(m,n)=(2,3)$ yields $\AD(G)\geq\frac{5}{3}$.
Otherwise, we must have $\cdset(G)=\{1,2\}$ and $\abs{G'}=2$, and following the steps in the proof of Lemma \ref{l:hammer} yields $\AD(G)=\frac{3}{2}$.
\end{proof}

\begin{rem}
In \cite{YC_minorant} the present author studied a generalization of $\AD(G)$ to the setting of virtually abelian groups, and showed that $\AD(G)=\frac{3}{2}$ if and only if $\gpind{G}{Z(G)}=4$.
The proof goes via a version of Corollary \ref{c:minimal AD}, but substantial work is required since $G$ may be infinite. It is therefore worth noting that when $G$ is finite, there is a much simpler proof of this equivalence; details are given in Appendix~\ref{app:minimal-AD}.
\end{rem}

We saw in the proof of Corollary \ref{c:minimal AD} that if $\AD(G)> \frac{3}{2}$, then either $\mindeg(G)\geq 3$ or $\abs{G'}\geq 3$.
The example of $S_3$ shows that we can have $\mindeg(G)=2$ and $\abs{G'}=3$. In contrast, the next result shows that we can never have $\mindeg(G)=3$ and $\abs{G'}=2$.

\begin{lem}\label{l:|G'|=2}
Let $G$ be a group with $\abs{G'}=2$. If $\varphi\in\Irr(G)$ and $\cd{\varphi}>1$, then $\cd{\varphi}$ is even.
\end{lem}

Lemma \ref{l:|G'|=2} follows from more precise results of Miller, stated in the introduction of~\cite{Miller_|G'|=2}. His presentation is rather terse and uses the finiteness of $G$ in an essential way. We provide a direct proof of Lemma \ref{l:|G'|=2} in Appendix~\ref{app:force even degree}, which also works for (irreducible, finite-dimensional, unitary) representations of infinite groups.

\begin{prop}\label{p:2 not in cdset}
Let $G$ be non-abelian. If $2\notin\cdset(G)$, then $\AD(G)\geq \frac{7}{3}$.
\end{prop}

\begin{proof}
We split into two cases.
If $\abs{G'}\geq 3$, then using Lemma \ref{l:hammer} with $m=2$ and $n=3$ gives $\AD(G)\geq \frac{7}{3}$.
If $\abs{G'}=2$, then $3\notin\cdset(G)$ by Lemma \ref{l:|G'|=2}, and so $\mindeg(G)\geq 4$; using Lemma \ref{l:hammer} with $m=4$ and $n=2$ gives
$\AD(G) \geq \frac{5}{2} > \frac{7}{3}$.
\end{proof}

In both cases of the proof, the lower bounds are sharp; see Example \ref{eg:extraspecial} below for details.
\end{section}

\begin{section}{The proof of Proposition \ref{p:new lower bound}}
\label{s:proof of main bound}
For a finite set $X$ and a function $f:X\to\Cplx$, we write $\supp(f)$ for the \dt{support of $f$}, that~is, the set $\{x\in X \colon f(x)\neq 0\}$.

\begin{lem}[The ``$\cL$-orbit method'' for lower bounds]\label{l:orbit method}
Let $\varphi\in\Irr(G)$ and let $n=\cd{\varphi}$. Let $K$ be the normal subgroup of $G$ generated by $\supp(\varphi)$. Then
\[
\abs{ \Irr_n(G) } \geq
\abs{\Irr_1(G)}\, \gpind{G}{K}^{-1} .
\]
\end{lem}

\begin{proof}
To simplify notation, let $\cL=\Irr_1(G)$.
$\cL$ is a group with respect to pointwise product, and multiplication of characters  defines a group action $\cL\times \Irr_n(G)\to\Irr_n(G)$ for each~$n$.
The $\cL$-orbit of $\varphi$ is a subset of $\Irr_n(G)$ and it has size $\abs{\cL}\, \abs{\Stab_{\cL}(\varphi)}^{-1}$.

Let $\bbT$ denote the set of complex numbers of unit modulus, viewed as a group with respect to multiplication.
Observe that each $\gamma\in\cL$ is $\bbT$-valued, and that
\[
\Stab_{\cL}(\varphi) = \{ \gamma\in\cL \colon \gamma\varphi=\varphi \}
= \{ \gamma\in\cL \colon \gamma(x)=1\text{ for all $x\in \supp(\varphi)$} \}\,,
\]
which is the set of group homomorphisms $G\to\bbT$ that factor through $G\to G/K$.
Therefore, writing $A$ for the abelianization of $G/K$, we have
\[
  \abs{\Stab_{\cL}(\varphi)} =\abs{A} \leq \abs{ G/K} = \gpind{G}{K}\,,
 \]
 and so the $\cL$-orbit of $\varphi$ has at least $\abs{\cL}\,\gpind{G}{K}^{-1}$ elements. The result now follows.
\end{proof}

\begin{cor}\label{c:crude lower bound}
Let $n\in\Nat$. If $\Irr_n(G)$ is non-empty, then $\abs{\Irr_n(G)}\geq n^{-2}\,\abs{\Irr_1(G)}$.
\end{cor}

\begin{proof}
Pick some $\varphi\in \Irr_n(G)$, and let $K$ be the normal subgroup of $G$ generated by $\supp(\varphi)$.
Since $\varphi$ is irreducible, $\ip{\varphi}{\varphi}=1$. Therefore, since $\abs{ \varphi(x)} \leq \cd{\varphi}=n$ for all $x\in G$,
\[
n^2 \abs{\supp(\varphi)} \geq \sum_{x\in G} \abs{\varphi(x)}^2 = \abs{G}\,.
\]
Hence $\gpind{G}{K} \leq \abs{G}\, \abs{\supp(\varphi)}^{-1} \leq n^2$. Applying Lemma \ref{l:orbit method}, the result follows.
\end{proof}

\begin{rem}\label{r:sharp}
Although the estimates in the proof of Lemma \ref{l:orbit method} are potentially wasteful, the resulting lower bound in Corollary \ref{c:crude lower bound} is sharp.
For if $G$ is an extraspecial group of order $2^{2k+1}$, it has exactly $2^{2k}$ characters of degree~$1$, and a single irreducible character of degree~$2^k$.
On the other hand, it will be important later that in certain situations we can do significantly better (Lemma \ref{l:probabilistic} below).
\end{rem}

\begin{prop}\label{p:maxdeg geq 4}
If $G$ is non-abelian, then $\AD(G) \geq 2+ (\maxdeg(G)-3)\abs{G'}^{-1}$.
\end{prop}

\begin{proof}
Let $d=\maxdeg(G)$.
Since $\abs{G} =\sum_{n=1}^d n^2 \abs{\Irr_n(G)}$,
\[ \begin{aligned}
\AD(G) -2
&= \frac{1}{\abs{G}}\sum_{n=1}^d (n^3-2n^2) \abs{\Irr_n(G)} \\
& \geq -\frac{ \abs{\Irr_1(G)} }{\abs{G}} +
(d^3-2d^2) \frac{\abs{\Irr_d(G)} }{\abs{G}} \,.
\end{aligned} \]
Since $\Irr_d(G)$ is non-empty, applying Corollary \ref{c:crude lower bound}  gives the desired inequality.
\end{proof}

The rest of this section deals with cases where $\cdset(G)=\{1,2,3\}$. We require a property of $\AD$ that is not obvious from the definition, but which seems to be crucial to understanding its behaviour.

\begin{prop}[Johnson]\label{p:monotone}
$\AD$ is monotone with respect to subgroup inclusion. That is, if $H\leq G$ then $\AD(H)\leq \AD(G)$.
\end{prop}

\begin{rem}
Proposition \ref{p:monotone} follows from results in  \cite[\S4]{BEJ_AG} concerning ``amenability constants'' of Fourier algebras, or from the general theory in \cite[\S2]{YC_minorant}.
One can give a direct proof, based on considering the induction of characters from $H$ to $G$, see the author's {\itshape MathOverflow} question \cite{MO_monotone} and the comments and answers. It is quite possible that a direct proof along these lines was already known to Johnson.
\end{rem}

\begin{prop}\label{p:index 2 inherits}
Let $G$ be a group such that $\AD(G)< \frac{7}{3}$.
Suppose that $\cdset(G)=\{1,2,3\}$; then every index-$2$ subgroup $H\leq G$ also satisfies $\cdset(H)=\{1,2,3\}$.
\end{prop}

The proof of Proposition \ref{p:index 2 inherits} requires some general facts about index-$2$ subgroups,
 which we state in a separate lemma for convenience.

\begin{lem}[Character degrees of index-2 subgroups]
\label{l:cd for index 2}
Let $H\leq G$ be an index-$2$ subgroup.
Then $\maxdeg(H)\leq\maxdeg(G)$ and $\cdset(G)\subseteq\cdset(H)\cup 2\cdset(H)$.
\end{lem}

Both parts of the lemma are standard results. For completeness, we quickly sketch their proofs.

\begin{proof}
Given $\psi\in\Irr(H)$, let $\varphi\in\Irr(G)$ be one of the irreducible summands of $\Ind^G_H\psi$.
By Frobenius reciprocity, $\psi$ is contained in $\Res{\varphi}{H}$, so $\cd{\psi}\leq \cd{\Res{\varphi}{H}}=\cd{\varphi}\leq\maxdeg(G)$.
This proves the first claim.

For the second claim, let $\varphi\in \Irr(G)$. If $\Res{\varphi}{H}$ is irreducible, then $\cd{\varphi}\in \cdset(H)$.
If not, then it follows from Clifford theory (or direct arguments using Frobenius reciprocity) that $\Res{\varphi}{H}$ splits as the sum of two irreducible characters, say $\beta_1$ and $\beta_2$, which satisfy $\varphi=\Ind^G_H\beta_1=\Ind^G_H\beta_2$.
In particular, $\cd{\varphi}=2\cd{\beta_1}\in 2\cdset(H)$.
\end{proof}

\begin{proof}[Proof of Proposition \ref{p:index 2 inherits}]
By monotonicity of $\AD$ (Proposition \ref{p:monotone}), $\AD(H)\leq \AD(G)< \frac{7}{3}$.
Hence by Lemma \ref{l:cd for index 2}, $\maxdeg(H)\leq 3$ and $3\in\cdset(H)$.
Since $H$ is non-abelian and $\AD(H)< \frac{7}{3}$, the contrapositive of Proposition \ref{p:2 not in cdset} implies that $2\in\cdset(H)$.
\end{proof}

We now observe that $2$-dimensional irreducible representations of $G$ can be used to produce $3$-dimensional representations with useful properties.
In what follows, $\veps$ denotes the constant function~$1$, regarded as the trivial representation of a group.

\begin{lem}\label{l:implies index 2 subgroup}
Let $\pi$ be a $2$-dimensional irreducible representation of $G$, and let $\pi^\ast$ denote its contragredient.
\begin{romnum}
\item\label{li:trivial summand}
$\veps$ occurs in $\pi\otimes\pi^\ast$ with multiplicity $1$.
\item\label{li:sign is summand}
Let $\rho$ be the summand in $\pi\otimes\pi^\ast$ complementary to $\veps$.
Suppose that $\rho$ is reducible. Then $G$ has an index-$2$ subgroup.
\end{romnum}
\end{lem}

This is surely not a new observation, but since we are unaware of a precise reference, a full proof is given below.

\begin{proof}
Part \ref{li:trivial summand} follows from Schur's lemma. (Alternatively, let $\psi=\Tr\pi$; then the multiplicity of $\veps$ in $\pi\otimes\pi^\ast$ is equal to
 $\ip{\psi\overline{\psi}}{\veps} = \abs{G}^{-1} \sum_{x\in G} \psi(x)\overline{\psi(x)} = 1$.)

For part \ref{li:sign is summand}, let $\varphi=\Tr\rho$;  by part~\ref{li:trivial summand}, $\varphi$ is real-valued and $\ip{\varphi}{\veps}=0$. We claim that there exists a real-valued character on~$G$ of degree~$1$ occuring as a summand of $\varphi$.
Assuming such a character exists, it may be viewed as a group homomorphism $\sigma:G\to\{\pm1\}$.
Since $\veps$ is not a summand of $\varphi$, we know that $\sigma\neq\veps$, and so $\ker\sigma$ is an index-$2$ subgroup of~$G$, as required.

To prove the claim, note that since $\rho$ has degree $3$ and is reducible, its decomposition into irreducible characters includes at least one $\gamma\in\Irr_1(G)$. If $\gamma$ is real-valued, we are done. If not, then $\overline{\gamma}\neq\gamma$ and $\ip{\varphi}{\overline{\gamma}}= \overline{\ip{\varphi}{\gamma}}\geq 1$. Hence $\gamma$ and $\overline{\gamma}$ occur in $\varphi$ with multiplicity~$1$, and $\varphi=\gamma+\overline{\gamma}+\sigma$ where $\sigma\in\Irr_1(G)$ is real-valued.
\end{proof}


If $G$ has no index-$2$ subgroups and $2\in\cdset(G)$,
 then by Lemma~\ref{l:implies index 2 subgroup}, for each $\psi\in\Irr_2(G)$ the character $\beta\defeq\psi\overline{\psi}-\veps$ must be irreducible; and because $\beta$ is a ``small perturbation'' of a non-negative character, we can obtain improved lower bounds on $\abs{\supp(\beta)}$,  allowing us to apply Lemma \ref{l:orbit method} more effectively. It turns out that the relevant estimates have nothing to do with group structure, so we present them as a separate lemma.

\begin{lem}\label{l:probabilistic}
Let $X$ be a finite non-empty set, and let $d\geq 1$. Suppose that $f:X\to [-1,d]$ has mean~$0$ and variance~$1$, i.e.
\[
\sum_{x\in X} f(x) = 0 \text{ and } \sum_{x\in X} f(x)^2 = \abs{X}\,.
\]
Then $\abs{\supp(f)}\geq d^{-1}\abs{X}$.
\end{lem}

\begin{proof}
Fix some ``threshold value'' $c \in [0,d]$, to be determined later, and partition $\supp(f)$ as $N\cup P \cup R$ where:
\begin{itemize}
\item $N \defeq \{ x\in X \colon -1\leq f(x) < 0 \}$;
\item $P \defeq \{ x\in X \colon 0 < f(x) \leq c \}$;
\item $R \defeq \{ x\in X \colon c < f \leq d \}$.
\end{itemize}
%
Then, since $\sum_{x\in\supp(f)} f(x)^2 =\abs{X}$,
\begin{equation}\label{eq:decompose and bound}
\begin{aligned}
\abs{X} & = \sum_{x\in N} f(x)^2 + \sum_{x\in P} f(x)^2 + \sum_{x\in R}
f(x)^2  \\
& \leq \sum_{x\in N} \abs{f(x)} + c \sum_{x\in P} f(x) +
d \sum_{x\in R} f(x) \,.
\end{aligned}
\tag{$*$}
\end{equation}
On the other hand, since $\sum_{x\in \supp(f)} f(x)=0$,
\[
\sum_{x\in P} f(x) = \sum_{x\in N} \abs{f(x)}- \sum_{x\in R} f(x) \,,
\]
and substituting this into \eqref{eq:decompose and bound} yields
\[ \begin{aligned}
\abs{X}  \leq (c+1) \sum_{x\in N} \abs{f(x)} + (d-c) \sum_{x\in R} f(x) 
& \leq (c+1) \abs{N}  + d(d-c) \abs{R}  \\
& \leq \max(c+1, d(d-c))\, \abs{\supp(f)}  \,.
\end{aligned}
\]
Taking $c=d-1$ gives $\abs{X}\leq d\abs{\supp(f)}$ as required.
\end{proof}

\begin{prop}\label{p:Irr3 better bound}
Suppose that $G$ has no index-$2$ subgroups but has an irreducible representation of degree~$2$.
 Then $\abs{\Irr_3(G)}\geq \frac{1}{3}\abs{\Irr_1(G)}$.
\end{prop}

\begin{proof}
Let $\psi\in\Irr_2(G)$ and let $\beta=\psi\overline{\psi}-\veps$. We observe that:
\begin{itemize}
\item $\beta$ takes values in $[-1,3]$, since $0\leq \abs{\psi(x)}^2\leq 4$ for all $x\in G$;
\item $\ip{\beta}{\veps}=0$, by Lemma \ref{l:implies index 2 subgroup}\ref{li:trivial summand};
\item $\ip{\beta}{\beta}=1$, since $\beta$ is irreducible by Lemma \ref{l:implies index 2 subgroup}\ref{li:sign is summand}.
\end{itemize}
Hence, by Lemma~\ref{l:probabilistic}, $\abs{\supp(\beta)}\geq \frac{1}{3}\abs{G}$, and applying Lemma~\ref{l:orbit method} completes the proof.
\end{proof}

\begin{rem}
In general, the  bound in Proposition \ref{p:Irr3 better bound} can not be improved. To see this, take $G=\SL(2,3)$. Then $\cdset(G)=\{1,2,3\}$ and $\abs{\Irr_1(G)}=3 = 3\abs{\Irr_3(G)}$, while $\gpind{G}{G'}=3$ (so that $G$ cannot quotient onto the two-element group).
\end{rem}

\begin{proof}[Proof of Proposition \ref{p:new lower bound}]
Let $G$ be a group with $\maxdeg(G)\geq 3$.
If $\maxdeg(G)\geq 4$, then $\AD(G)\geq 2 + \abs{G'}^{-1}$ by Proposition~\ref{p:maxdeg geq 4}. So we assume henceforth that $\maxdeg(G)=3$. Note that this implies $\abs{G'}\geq 3$, by Lemma~\ref{l:|G'|=2}.
Moreover, if $\cdset(G)=\{1,3\}$, then by Proposition~\ref{p:2 not in cdset}, $\AD(G) \geq \frac{7}{3} \geq 2+ \abs{G'}^{-1}$.

It only remains to deal with the cases where $\cdset(G)=\{1,2,3\}$. If $\AD(G)\geq \frac{7}{3}$ then we are done, as before. So we may assume that $\cdset(G)=\{1,2,3\}$ and $\AD(G)< \frac{7}{3}$. Put $H_0=G$ and apply the following recursive procedure:
 if $n\in\Nat$ and $H_{n-1}$ has an index-$2$ subgroup, choose $H_n$ to be such a subgroup; otherwise, stop.
 Note that at each stage, Proposition \ref{p:index 2 inherits} ensures that $\cdset(H_n)=\{1,2,3\}$.

Since $G$ is finite this procedure must terminate; let $H$ be the last subgroup in this sequence. Since $\cdset(H)=\{1,2,3\}$,
\[ \begin{aligned}
\AD(H) & = \frac{1}{\abs{H}} \Bigl( \abs{\Irr_1(H)} + 8 \abs{\Irr_2(H)} + 27 \abs{\Irr_3(H)} \Bigr), \text{ and} \\
1 &  = \frac{1}{\abs{H}} \Bigl( \abs{\Irr_1(H)} + 4 \abs{\Irr_2(H)} + 9 \abs{\Irr_3(H)} \Bigr) \,.
\end{aligned} \]
Hence $\AD(H) = 2-\abs{H}^{-1} \abs{\Irr_1(H)} + 9 \abs{H}^{-1} \abs{\Irr_3(H)}$. On the other hand, since $H$ has no index-$2$ subgroups, it satisfies the hypotheses of Proposition \ref{p:Irr3 better bound}, and so
\[
\AD(H)  \geq 2 + \frac{2\abs{\Irr_1(H)}}{\abs{H}} = 2+ \frac{2}{\abs{H'}}\,.
\]
As $\AD(G)\geq\AD(H)$ (Proposition \ref{p:monotone}) and $\abs{G'}\geq \abs{H'}$, we conclude that $\AD(G) \geq 2 + 2\abs{G'}^{-1}$, which completes the proof of Proposition \ref{p:new lower bound}.
\end{proof}

\end{section}

\begin{section}{Further examples and implications of small values}

\begin{subsection}{Values of $\AD$ for particular groups}
\label{ss:particular}
We present three families of groups with rather different properties (nilpotent, solvable with trivial centre, and quasi-simple), where one obtains rather simple formulas for the $\AD$-constants in each family. In each case, the ratio $\AD(G) \maxdeg(G)^{-1}$ converges to $1$ as $\abs{G}\to\infty$. 

\begin{eg}[Extraspecial $p$-groups]\label{eg:extraspecial}
Let $p$ be a prime and let $n\in\Nat$. If $G$ is an extraspecial $p$-group of order $p^{2n+1}$, then the degrees of its irreducible characters and their multiplicities are well-documented.
Namely, $G$ has exactly $p^{2n}$ characters of degree~$1$, and exactly $p-1$ irreducible characters of degree~$p^n$. Hence
\begin{equation}
  \AD(G)= \frac{p^{2n}\cdot 1^3  + (p-1)p^{3n}}{p^{2n+1}} = p^{n-1}(p-1) + \frac{1}{p} \,.
\end{equation}
We note two particular cases, relevant to Proposition \ref{p:2 not in cdset}.
If $p=2$ and $n=2$, then $\cdset(G)=\{1,4\}$ and $\AD(G)=\frac{5}{2}$.
If $p=3$ and $n=1$, then $\cdset(G)=\{1,3\}$ and $\AD(G)=\frac{7}{3}$.
\end{eg}

\begin{eg}[Affine groups of finite fields]\label{eg:AD of Aff(F)}
For $q$ a prime power $\geq 3$, let $\bbF_q$ denote the finite field with $q$ elements, and consider the natural semidirect product $\bbF_q\rtimes \bbF_q^\times$ (sometimes referred to as the affine group or ``$ax+b$ group'' of $\bbF_q$). This group has exactly $q-1$ characters of degree~$1$, and a single irreducible character of degree~$q-1$. Hence
\begin{equation}
\AD(\bbF_q\rtimes\bbF_q^\times) = \frac{(q-1)\cdot 1^3 + (q-1)^3}{q(q-1)} = q- 2+ \frac{2}{q}\,.
\end{equation}
Note that when $q=3$ this group is isomorphic to the dihedral group of order $6$, and its $\AD$-constant is $\frac{5}{3}$; this matches the calculation in Example \ref{eg:dihedral}.
\end{eg}
\begin{eg}[Special linear groups of degree $2$]
Let $q$ be a prime power, and let $\SL(2,q)$ denote the special linear group of degree~$2$ over the finite field with $q$ elements; this has order $q^3-q$.

For $q$ even, put $q=2r$; then $\Irr(\SL(2,q))$ is the union of four pairwise disjoint sets $X_1$, $X_{q-1}$, $X_q$ and $X_{q+1}$,
 where each member of $X_j$ has degree~$j$, and
\[
\hbox{$\abs{X_1}=1$;
$\abs{X_{2r-1}}=r$;
$\abs{X_{2r}}=1$;
$\abs{X_{2r+1}}=r-1$.}
\]

For $q$ odd, put $q=2r+1$; then $\Irr(\SL(2,q))$ is the union of six pairwise disjoint sets $X_1$, $X_r$, $X_{r+1}$, $X_q$ and $X_{q+1}$,
 where each member of $X_j$ has degree~$j$, and
\[
\hbox{$\abs{X_1}=1$;
$\abs{X_r} = 2$;
$\abs{X_{r+1}} = 2$;
$\abs{X_{2r}} = r$;
$\abs{X_{2r+1}} = 1$;
$\abs{X_{2r+2}} = r-1$.}
\]

By brute-force calculation, we eventually obtain
\begin{equation}
  \AD(\SL(2,q))  = \left\{
    \begin{aligned}
    \frac{q^3-3}{q^2-1}  & = q - \frac{1}{q-1} + \frac{2}{q+1} & \text{for $q$ even,} \\
    \frac{2q^3-q^2-9}{2(q^2-1)}  & = q - \frac{1}{2} -	\frac{2}{q-1} + \frac{3}{q+1}&\text{for $q$ odd.}
  \end{aligned}
\right.
\end{equation}

\end{eg}

The next set of examples was suggested to the author by P. Levy.

\begin{eg}[Finite subgroups of $\SO(3)$ and $\SU(2)$]
\label{eg:AD double covers}
We ignore the cyclic groups and dihedral groups, and their double covers inside $\SU(2)$, since these are covered by previous results. So there are only three new examples to consider.
In the following list, when we refer to the ``character degrees'' of a group $H$, we mean ``the degrees of its irreducible characters, listed with multiplicity''.

\begin{alphnum}
\item\label{li:bin-tet}
The alternating group $A_4$ has character degrees $1,1,1,3$.
Its double cover is the binary tetrahedral group $2T\cong \SL(2,3)$, whose character degrees are $1, 1, 1, 2, 2, 2, 3$. Thus
\[
\AD(A_4)= \frac{30}{12} = \frac{5}{2}
\text{ and }
\AD(2T)= \frac{54}{24} =\frac{9}{4} < \AD(A_4).
\]

\item\label{li:bin-oct}
The symmetric group $S_4$ has character degrees $1,1,2,3,3$.
Its double cover is the binary octahedral group $2O$, whose character degrees are $1, 1, 2, 2, 2, 3, 3, 4$. Thus
\[
\AD(S_4) = \frac{64}{24} = \frac{8}{3}
\text{ and }
\AD(2O) = \frac{144}{48} = 3 > \AD(S_4).
\]

\item\label{li:bin-icos}
The alternating group $A_5$ has character degrees $1, 3, 3, 4, 5$.
Its double cover is the binary icosahedral group $2I\cong\SL(2,5)$, whose character degrees are $1, 2, 2, 3, 3, 4, 4, 5, 6$. Thus
\[
\AD(A_5) = \frac{244}{60} = \frac{61}{15}
\text{ and }
\AD(2I) = \frac{540}{120} = \frac{9}{2} > \AD(A_5).
\]
\end{alphnum}

\end{eg}

\begin{rem}\label{r:AD bad quotient}
It was already known that although $\AD$ cannot increase when passing to subgroups, it can increase when passing to a quotient. For instance, in a ``note added in proof'' in \cite{LLW_SM96}, it is observed that the Schur cover of $A_6$ has an $\AD$ constant strictly smaller than that of the triple cover of $A_6$. However, Example \ref{eg:AD double covers}\ref{li:bin-tet} shows that they missed a much smaller example.
\end{rem}

\end{subsection}

\begin{subsection}{Structural consequences for $G$ of upper bounds on $\AD$}
\label{ss:miscell}

\begin{prop}[A cheap lower bound for $p$-groups]
Let $p$ be a prime. If $G$ is a non-abelian $p$-group, then $\AD(G) \geq p - 1 + \frac{1}{p}$. Equality is attained by an extraspecial $p$-group of order $p^3$.
\end{prop}

\begin{proof}
Since $G$ is a $p$-group, both $\mindeg(G)$ and $\abs{G'}$ are powers of $p$. Therefore both are $\geq p$, since $G$ is non-abelian. The rest follows from Lemma \ref{l:hammer} and the calculation in Example \ref{eg:extraspecial}.
\end{proof}

A similar idea can be used to control (sub)groups of odd order whose $\AD$-constant is small. The next result is a slightly stronger version of an observation by G. Robinson (personal communication).

\begin{lem}\label{l:odd order subgroup}
If $\AD(G)< \frac{7}{3}$ then every odd order subgroup of $G$ is abelian.
\end{lem}

\begin{proof}
We prove the contrapositive. Suppose that $G$ has a non-abelian subgroup $H$ that has odd order.
Since $\cd{\varphi}$ divides $\abs{H}$ for each $\varphi\in\Irr(H)$, we have $\mindeg(H)\geq 3$; since $\abs{H'}$ divides $\abs{H}$, we have $\abs{H'}\geq 3$. Therefore, by monotonicity (Proposition \ref{p:monotone}) and Lemma~\ref{l:hammer},
\[ \AD(G) \geq \AD(H) \geq 1 + (3-1) \frac{2}{3} = \frac{7}{3} \,, \]
as required.
\end{proof}

\begin{cor}
If $G$ is nilpotent and $\AD(G)< \frac{7}{3}$, then $G$ is the product of a $2$-group and an abelian group of odd order.
\end{cor}

\begin{proof}
If $p$ is an odd prime, then by Lemma \ref{l:odd order subgroup} each $p$-Sylow subgroup of $G$ is abelian. But since $G$ is finite and nilpotent, it factorizes as the direct product of its Sylow subgroups.
\end{proof}

\begin{rem}
We can show by relatively elementary arguments that $\AD(G)> 4$ whenever $G$ is non-abelian and simple; since $\AD(A_5)=\frac{61}{15}$, this is already quite close to the optimal result. Although we shall obtain a stronger result in Section~\ref{s:nonsolvable}, we include the proof of the weaker result here as an illustration of our earlier method.

The main idea is similar to the proof of Lemma \ref{l:hammer}. Let $m=\mindeg(G)$. Since $m$ divides $\abs{G}=\sum_{\varphi\in\Irr(G)} \cd{\varphi}^2$ and since $\abs{\Irr_1(G)}=1$,
\[
\sum_{\varphi \in\Irr(G), \cd{\varphi}>m} \cd{\varphi}^2 \equiv -1\bmod m \,.
\]
Hence, there is at least one $\sigma\in\Irr(G)$ with $\cd{\sigma}\geq m+1$. Therefore
\[ \begin{aligned}
    \AD(G)- m & = \sum_{\varphi\in\Irr(G)} \frac{ (\cd{\varphi}-m) \cd{\varphi}^2 }{ \abs{G} } \\
    & \geq - \frac{m-1}{\abs{G}} + \frac{(\cd{\sigma}-m)\cd{\sigma}^2}{\abs{G}} \geq   - \frac{m-1}{\abs{G}} + \frac{(m+1)^2}{\abs{G}} > 0\,.
  \end{aligned}
\]
If $m\geq 4$ this immediately gives $\AD(G)>4$. So it only remains to deal with cases where $m=3$.
The simple finite subgroups of $\PGL(3,\Cplx)$ were determined by Blichfeldt in \cite{Blichfeldt_PGL3},
 and using his classification one can show that the only simple groups with $m=3$ are $A_5$ and $\PSL(2,7)$ (some further explanation is given in Appendix~\ref{app:bypass CFSG}). We saw in Example \ref{eg:AD double covers}\ref{li:bin-icos} that $\AD(A_5)=\frac{61}{15} > 4$, and since $\PSL(2,7)$ has character degrees $1,3,3,6,7,8$ we find that $\AD(\PSL(2,7))=\frac{563}{84}> 6$.
\end{rem}
\end{subsection}

\end{section}

\begin{section}{A sharp lower bound on $\AD$ for non-solvable groups}
\label{s:nonsolvable}
Our aim in this section is to prove Theorem \ref{t:implies solvable}: if $G$ is non-solvable, then $\AD(G) \geq \frac{61}{15}$. 
One difficulty if we rely on our existing tools is that although $\AD$ behaves well with respect to taking subgroups, it does not behave well with respect to taking quotients (c.f.~Remark~\ref{r:AD bad quotient}).

Instead, our proof is inspired by techniques used in \cite{TongViet_large} to prove an analogous ``threshold'' result for the quantity
\[
f(G) \defeq \frac{1}{\abs{G}} \sum_{\varphi\in\Irr(G)} \cd{\varphi}\,.
\]
The key in \cite{TongViet_large} is to exploit the inequality $f(G)^2\leq \cp(G)$, where the \dt{commuting probability} $\cp(G)$ is equal to $\abs{G}^{-1} \abs{\Irr(G)}$. (Strictly speaking, this is not the original definition of $\cp$, but the equivalence with this definition is well known.) The inequality relating $f$ with $\cp$ is immediate from Cauchy--Schwarz; in our setting, we can use H\"older's inequality to obtain an analogous relationship between $\AD$ and $\cp$, but in the opposite direction.

\begin{prop}\label{p:AD and cp}
For every $G$ we have $1\leq \AD(G)^2 \cp(G)$. Equality is strict if $G$ is non-abelian.
\end{prop}

\begin{proof}
Applying H\"older's inequality with conjugate exponents $\frac{3}{2}$ and $3$ gives
\[
\abs{G} = \sum_{\varphi\in\Irr(G)} \cd{\varphi}^2\cdot 1 \leq \left(\sum_{\varphi\in\Irr(G)} \cd{\varphi}^3 \right)^{\frac{2}{3}} \left( \sum_{\varphi\in\Irr(G)} 1^3\right)^{\frac{1}{3}} \;,
\]
and the inequality is strict unless every $\varphi\in\Irr(G)$ has the same degree, i.e.\ unless $G$ is abelian. The result now follows by dividing both sides by $\abs{G}$ and then cubing.
\end{proof}

\begin{rem}\label{r:hat-tip_GR}
The invariant $\cp$ has been intensively studied, and in particular it is shown in \cite[Theorem 11]{GR_cp} that if $\frac{1}{12} > \cp(G) >\frac{3}{40}$ then $G$ is solvable.
Although this result itself is not strong enough to imply Theorem \ref{t:implies solvable}, the ideas in its proof can be seen (refracted through the prism of \cite{TongViet_large}) in what follows.
\end{rem}

 The key advantage of working with $\cp$, compared with either $f$ or $\AD$, is that it behaves well with respect to both taking subgroups and taking quotients. In particular, we make crucial use of the following result.

\begin{lem}[Gallagher, \cite{Gallagher_cp}]
\label{l:gallagher}
Suppose that $N\unlhd G$. Then
\[
\min(\cp(G/N),\cp(N)) \geq \cp(G/N)\cp(N) \geq \cp(G) \,.
\]
\end{lem}

 It was observed by Dixon that the largest value of $\cp$ on simple non-abelian groups is attained at~$A_5$. We need some information about the second largest value attained by $\cp$ on this class of groups.
The following is a slightly stronger version of \cite[Lemma 2.3]{TongViet_large}

\begin{lem}\label{l:annoying}
There exists $\frac{1}{28} \leq \delta_0 \leq \frac{1}{20}$ with the following property: if $S$ is a finite non-abelian simple group and $\cp(S) > \delta_0$, then $S\cong A_5$ (and $\cp(S)=\frac{1}{12}$).
\end{lem}

By consulting the classification of finite simple groups (CFSG) and considering the minimal degrees of non-trivial irreducible characters, it can be shown that one can take $\delta_0=\frac{1}{28}$. (This is best possible since $\cp(\PSL(2,7))=\frac{1}{28}$.) 
We can show without resorting to the full CFSG that $\delta_0 = \frac{1}{20}$ works. Details are given in Appendix~\ref{app:bypass CFSG}: our approach invokes parts of the classification of finite subgroups of $\PGL(3,\Cplx)$ and $\PGL(4,\Cplx)$, given by Blichfeldt in the 1900s \cite{Blichfeldt_PGL4,Blichfeldt_PGL3}.

The following lemma is presumably standard knowledge, but it seems quicker to give an explanation than to look up a reference.

\begin{lem}\label{l:orbits of perfect}
Let $H$ be a perfect group.
\begin{alphnum}
\item\label{li:perfect to solvable} If $W$ is a solvable group and $\theta:H\to W$ is a homomorphism, then $\theta(H)=\{\gpid[W]\}$.
\item\label{li:orbits} If $X$ is any set on which $H$ acts, then each $H$-orbit in $X$ either has size $1$ or size $\geq 5$.
\end{alphnum}
\end{lem}

\begin{proof}
Part \ref{li:perfect to solvable} follows by induction on the derived series of~$W$. For part \ref{li:orbits}, observe that an $H$-orbit of size $n$ defines a homomorphism $\alpha:H \to S_n$ whose image acts transitively on the original orbit. If $n\leq 4$, then $S_n$ is solvable and so $\alpha(H)$ is trivial by part \ref{li:perfect to solvable}; this is only possible if~$n=1$.
\end{proof}

We now turn to the proof that the only perfect groups with commuting probability greater than $\frac{1}{20}$ are $A_5$ and $\SL(2,5)$ (Proposition \ref{p:sharper TV}). Our argument is patterned on the proof of  \cite[Lemma~2.4]{TongViet_large}, but since we need better bounds than those in Tong-Viet's paper, we take the opportunity to make some simplifications and give a more streamlined approach.

\begin{proof}[Proof of Proposition \ref{p:sharper TV}]
Let $S$ be the quotient of $H$ by any maximal proper normal subgroup. Then $S$ is simple (by maximality), non-trivial (by properness), and non-abelian (since $H$ is perfect). By Lemma~\ref{l:gallagher} $\cp(S)\geq \cp(H) > \frac{1}{20}$, so by Lemma \ref{l:annoying}, $S\cong A_5$. 

Thus we have a surjective homomorphism $H \to A_5$, with kernel $N$ say. If $N$ is trivial, there is nothing to prove;
so henceforth, we assume $\abs{N}\geq 2$, and aim to prove that $H\cong\SL(2,5)$.

By definition, $H$ is an extension of $A_5$ by the group $N$.
Suppose that we can show it is a \emph{central} extension; then
since $H$ is perfect, it must be a quotient of the Schur cover of $A_5$, which is isomorphic to $\SL(2,5)$.
Since $\abs{\SL(2,5)}=2\abs{A_5}\leq\abs{H}$, the quotient map from $\SL(2,5)$ onto $H$ must be injective, and we are done.

Therefore, it suffices to prove that $N\subseteq Z(H)$.
Let $k_H(N)$ denote the number of $H$-conjugacy classes contained in $N$.
As in the proof of \cite[Lemma 2.4]{TongViet_large}, we have the inequality

\begin{equation}\label{eq:needs machinery}
12 \cp(H) \leq \frac{k_H(N)}{\abs{N}}
\end{equation}

(We briefly sketch how this works. If $M$ is any finite group and $N\unlhd M$, then
 \cite[Lemma~1(iii)]{GR_cp}, which is actually proved in \cite[Remark A2']{KovRob}, tells us that
\[
\abs{\Irr(M)} \leq k_M(N) \sup_{B\leq M/N} \abs{\Irr(B)} \,.
\]
We then apply this with $M=H$, noting that $\gpind{H}{N}=60$, and appeal to the fact that each subgroup of $A_5$ has at most five distinct irreducible characters.)

Since $\cp(H)>\frac{1}{20}$, it follows from \eqref{eq:needs machinery} that $k_H(N) > \frac{3}{5}\abs{N}$. On the other hand, $k_H(N)$ counts the number of orbits for the conjugation action of $H$ on $N$. By Lemma \ref{l:orbits of perfect}, the size of each non-singleton orbit is at least $5$.
Therefore, if $F$ denotes the set of fixed points of the action,
\[
\abs{N} \geq 5(k_H(N)-\abs{F}) + \abs{F} = 5k_H(N)-4\abs{F}\,,
\]
and combining this with the previous lower bound on $k_H(N)$ gives
\[
\abs{F} \geq \frac{1}{4}\left( 5k_H(N) - \abs{N}\right) > \frac{1}{2}\abs{N}\,.
\]
Now observe that $F=Z(H)\cap N$. So by the previous inequality $\gpind{N}{Z(H)\cap N}<2$, which is only possible if $Z(H)\cap N =N$, and this completes the proof.
\end{proof}

From this, we can show that on the class of finite perfect groups, the $\AD$-constant is minimized at $A_5$. In fact, a more precise statement can be made.

\begin{cor}\label{c:smallest AD perfect}
Let $H$ be a non-trivial perfect group which satisfies $\AD(H)\leq 2\sqrt{5} \approxeq 4.47$. Then $H\cong A_5$ and $\AD(H)=\frac{61}{15}\approxeq 4.07$.
\end{cor}

\begin{proof}
By Proposition \ref{p:AD and cp}, $\cp(H) > \AD(H)^{-2} \geq \frac{1}{20}$. Hence, by Proposition \ref{p:sharper TV}, $H$ is isomorphic to either $A_5$ or $\SL(2,5)$. But the second possibility is excluded, since we saw in Example \ref{eg:AD double covers}\ref{li:bin-icos} that $\AD(\SL(2,5))=\frac{9}{2} > 2\sqrt{5}$.
\end{proof}

\begin{proof}[Proof of Theorem \ref{t:implies solvable}]
Since $G$ is not solvable, its derived series stabilizes at some subgroup $H\leq G$ that is perfect and non-trivial. By monotonicity, $\AD(G)\geq\AD(H)$; and by Corollary \ref{c:smallest AD perfect}, we have $\AD(H)  \geq \min\left(2\sqrt{5}, \frac{61}{15}\right)=\frac{61}{15}$.
\end{proof}

\end{section}

%
%
%

\subsection*{Acknowledgements}
The results presented here are part of a wider programme to understand the amenability constants of Fourier algebras of locally compact groups, some of which was presented at the 2022 meeting of the \textit{Canadian Abstract Harmonic Analysis Symposium}. I~would like to thank the organizers and participants of that meeting for their interest; particular thanks are due to Brian Forrest and John Sawatzky.

Closer to home, thanks are due to Paul Levy for several helpful conversations. I~also thank Geoff Robinson for his interest in the topic of this paper at an earlier stage, and for various messages of encouragement.
The paper has also benefited from a close reading by an anonymous referee, whose suggestions have improved the clarity and consistency of the presentation in several places.

Finally, I would like to acknowledge the work of the contributors and maintainers of the {\itshape GroupPropsWiki}, and all those on {\itshape MathOverflow} who have been willing to share ideas and expertise with a Bear of Very Little Brain.

\appendix
\newpage
\begin{section}{Easier proofs of some known results}

\begin{subsection}{Finite groups with two character degrees and derived subgroup of order $2$}
\label{app:minimal-AD}
The groups described in the title are, by Corollary \ref{c:minimal AD}, those finite non-abelian groups where $\AD$ attains its minimum value. In this section, we give a quick proof that these groups are precisely those in which the centre has index~$4$.

Let $G$ be a finite group with $\cdset(G)=\{1,2\}$ and $\abs{G'}=2$.
We have $\Irr(G)=\Irr_1(G)\cup\Irr_2(G)$; let $l=\abs{\Irr_1(G)}$ and $m=\abs{\Irr_2(G)}$.
Also, since every conjugacy class injects into $G'$ and $\abs{G'}=2$, each conjugacy class in $G$ has size $1$ or $2$. Let $s=\abs{Z(G)}$ and let $n$ be the number of conjugacy classes of size~$2$.

Note that $\abs{G}=\sum_{\varphi\in\Irr(G)}\cd{\varphi}^2 =l+4m$. Since $\abs{G'}=2$, we have $\abs{G}=2l$, and so $l=4m$.
On the other hand, $\abs{G}= s+2n$, while $s+n=l+m$ since the character table is square.  Therefore,
\[ s+ 2n = 8m \quad,\quad s+n = 5m. \]
Solving for $s$ and $n$ yields $s=2m$ and $n=3m$. In particular, we conclude that
\[
\gpind{G}{Z(G)} = \frac{8m}{2m} = 4\,.
\]

Conversely, suppose that $\gpind{G}{Z(G)}=4$. The argument that follows is essentially the same as in \cite{YC_minorant}, but we include the details for sake of completeness.

Note that $G/Z(G)$ cannot be cyclic (otherwise by lifting the generator we would find that $G$ is abelian) and hence it is isomorphic to $C_2\times C_2$. Pick two generators for $G/Z(G)$ and lift them to $x,y\in G$. Then $x^2$, $y^2$ and $[x,y]$ all belong to $Z(G)$.
Since $G= Z(G) \cup xZ(G) \cup yZ(G) \cup xy Z(G)$, a short case-by-case analysis shows that every commutator in $G$ equals either $\gpid$ or $[x,y]$. In particular $\abs{G'}=2$.

Moreover, $A=Z(G)\cup x Z(G)$ is an Abelian subgroup of $G$ with index $2$. Hence, as shown in the proof of Lemma \ref{l:cd for index 2}, $\cdset(G)\subseteq \{1,2\}$. Since $G$ is non-abelian, this inclusion of sets is an equality.
\end{subsection}

\begin{subsection}{A self-contained proof of Lemma \ref{l:|G'|=2}}
\label{app:force even degree}
We give a proof of Lemma \ref{l:|G'|=2}, which works even for infinite groups.
Thus, for this subsection only, we let $G$ be a not-necessarily-finite group, and we suppose that $\abs{G'}=2$.

Let $\pi$ be a finite-dimensional, unitary, irreducible representation of $G$ with dimension $d\geq 2$. (When $G$ is finite, every irreducible representation of $G$ is automatically finite-dimensional and is equivalent to a unitary representation.)
Our aim is to show that $d$ is even.

\begin{lem}\label{l:G' is central}
Let $G'=\{\gpid,z\}$. Then $z\in Z(G)$.
\end{lem}

\begin{proof}
If $\alpha\in \operatorname{Aut}(G)$, then $\alpha(G')=G'$ and $\alpha'(\gpid)=\gpid$, hence $\alpha(z)=z$. Now take $\alpha$ to be an arbitrary inner automorphism of $G$.)
\end{proof}

\begin{lem}\label{l:squares are central}
Let $g\in G$. Then $g^2\in Z(G)$.
\end{lem}

\begin{proof}
Let $g,x\in G$. Then $gxg^{-1}=[g,x]x$. Since $[g,x]$ is central (by Lemma \ref{l:G' is central}) and $[g,x]^2=\gpid$, 
\[
g^2xg^{-2} = g \bigl([g,x]x\bigr)g^{-1} = [g,x](gxg^{-1})=[g,x][g,x] x = x\,.
\]
Thus $g^2$ is central.
\end{proof}

\begin{proof}[Proof that $d$ is even]
Recall that every non-trivial commutator in $G$ is equal to $z$.
Since $\pi$ is irreducible and $d\geq 2$, $\pi(G)$ is not Abelian, hence $\pi(z)\neq I_\pi$.
By Lemma \ref{l:G' is central} and Schur's lemma, $\pi(z)$ is a scalar multiple of $I_\pi$; since $z^2=\gpid$, it follows that $\pi(z)=-I_\pi$.

Since $G$ is non-abelian, there exist $x,y\in G$ that do not commute. Since $xyx^{-1}=zy$, we have $\pi(x)\pi(y)\pi(x)^{-1}=-\pi(y)$.
On the other hand, by Lemma \ref{l:squares are central} and Schur's lemma, $\pi(y^2)$ is a scalar multiple of $I_\pi$. Pick $\lambda\in \bbT$ such that $\pi(y)^2=\lambda^2 I_\pi$; then $U\defeq\lambda^{-1}\pi(y)$ is an involution in $\operatorname{Lin}(H_\pi)$ and $U$ is conjugate to $-U$.

Since $U$ is an involution it has exactly $d$ eigenvalues counted with multiplicity, and these eigenvalues belong to $\{-1,1\}$; but since $U$ is conjugate to $-U$, the eigenvalues $-1$ and $1$ must occur with equal multiplicity, $m$ say. Thus $d=2m$. (Alternatively, observe that $\frac{1}{2}(I_d+U)$ is an idempotent which has trace equal to $\frac{d}{2}$, which again forces $d$ to be even.)
\end{proof}

\end{subsection}

\begin{subsection}{A proof of Lemma \ref{l:annoying} with $\delta_0= 1/20$}
\label{app:bypass CFSG}
We follow the general strategy seen in the proofs of \cite[Lemma 2.3]{TongViet_large} and \cite[Theorem~11]{GR_cp}.
Suppose that $S$ is simple and non-abelian. Writing $m$ for $\mindeg(S)$, we have
\[
\abs{S} -1 \geq \sum_{\varphi\in\Irr(G), \cd{\varphi}>1} \cd{\varphi}^2 \geq \left(\abs{\Irr(S)}-1\right)m^2 = m^2\cp(S)\abs{S}- m^2\,,
\]
and rearranging gives
$m^2-1 \geq (m^2\cp(S)-1) \abs{S}$. If we are given an explicit $\delta_0>0$ such that $\cp(S)>\delta_0$, it follows that
\begin{equation}\label{eq:cheap leverage}
\frac{1}{\abs{S}}
 \geq \frac{m^2\cp(S)-1}{m^2-1}
 > \frac{m^2\delta_0-1}{m^2-1}\,.
\end{equation}
Provided that $m^2\delta_0 > 1$, the inequality \eqref{eq:cheap leverage} gives an explicit upper bound on $\abs{S}$.

Thus, in cases where $m$ is sufficiently large, $S$ belongs to some small list of known examples, and in each case we can see by inspection that $m$ is actually small (giving a contradiction). Separate {\it ad~hoc} arguments are then needed to deal with the cases where $m$ is ``small''.

In \cite[Lemma 2,3]{TongViet_large}, this strategy is used with $\delta_0 = \frac{16}{225}$, and so the easy part of the argument works for all $m\geq 4$; the only remaining cases are those with $m=3$, and these are covered by the following result.

\begin{thm}[Blichfeldt, implicitly]
\label{t:Blichfeldt 3-dim}
Let $S$ be a finite simple group with an irreducible representation of degree~$3$. Then $S\cong A_5$ or $S\cong \PSL(2,7)$.
\end{thm}

Inspecting the proof of \cite[Lemma 2.3]{TongViet_large}, the ``threshold value'' stated there can be improved from~$\frac{16}{225}$ to  $\frac{1}{15}$, provided that we know the simple groups of order $\leq 225$. However, the methods in that paper cannot reach $\frac{1}{16}$ (since we require $m^2\delta_0>1$), and for our eventual application to Proposition \ref{p:sharper TV}, we require $\delta_0\leq\frac{1}{20}$. We therefore need the following additional result.

\begin{thm}[Blichfeldt, implicitly]
\label{t:Blichfeldt 4-dim}
Let $S$ be a finite simple group with an irreducible representation of degree~$4$. Then $S\cong A_5$.
\end{thm}

For the reader who wishes to consult the original papers, we provide some details of how the theorems stated above are derived from the results stated in \cite{Blichfeldt_PGL4,Blichfeldt_PGL3}.

\begin{proof}[Proofs of Theorem \ref{t:Blichfeldt 3-dim} and \ref{t:Blichfeldt 4-dim}]
Let $d\in\{3,4\}$, and let $S$ be a finite simple group with an irreducible representation of degree~$d$. Then the image of $S$ under this representation can be identified with a finite subgroup $\widetilde{S}\leq SU(d)$ that acts irreducibly on $\Cplx^d$. Let $S_0$ be the image of $\widetilde{S}$ in $\PGL_d(\Cplx)$. Of course, $S\cong \widetilde{S}\cong S_0$.

In  the language of \cite[p.~553]{Blichfeldt_PGL3} and \cite[p.~205]{Blichfeldt_PGL4}, $S_0$ is primitive: this follows from the fact that simple groups cannot act non-trivially on sets of size $\leq 4$, c.f.\ Lemma \ref{l:orbits of perfect}.

\paragraph{The case $d=3$.}
The primitive simple finite subgroups of $\PGL_3(\Cplx)$ are determined up to isomorphism in \cite[Section 24]{Blichfeldt_PGL3} (relying on previous work of Maschke): any such subgroup must be isomorphic to $A_5$, $\PSL(2,7)$ or $A_6$. Moreover, it is observed that in the last case $A_6$ cannot be lifted from $\PGL_3(\Cplx)$ up to $\GL_3(\Cplx)$; thus $\widetilde{S}$ must be isomorphic to either $A_5$ or $\PSL(2,7)$, and this completes the proof of Theorem~\ref{t:Blichfeldt 3-dim}. 

\paragraph{The case $d=4$.}
The primitive simple finite subgroups of $\PGL_4(\Cplx)$ are determined up to isomorphism in \cite[Section III]{Blichfeldt_PGL4}; the list appears as items $22^\circ$--$27^\circ$  on pages 225--226 of that paper, and consists (in modern notation) of:
$A_5$, $A_6$, $A_7$, $\PSL(2,7)$, and $\PSp(4,3)$.
Blichfeldt does not state for which of these $S_0$ the corresponding ``lift'' in $\GL(4,\Cplx)$ is simple, but if we invoke known character tables for these groups then we see that the only possibility  for $\widetilde{S}$ is $A_5$ (none of the others have irreducible representations of degree~$4$), and this completes the proof of Theorem~\ref{t:Blichfeldt 4-dim}.
\end{proof}

\begin{proof}[Proof of Lemma \ref{l:annoying} with $\delta_0=\frac{1}{20}$.]
Let $S$ be non-abelian and simple, and let $m=\mindeg(S)$. Suppose that $\cp(S)>\frac{1}{20}$.
We start by showing that this forces $m\leq 4$. For, if $m\geq 5$, taking $\delta_0=\frac{1}{20}$ in \eqref{eq:cheap leverage} gives
\[
\frac{1}{\abs{S}}
 > \left(\frac{m^2}{20}-1 \right) \frac{1}{m^2-1}
\geq
 \left(\frac{25}{20}-1\right)\frac{1}{24} = \frac{1}{96}\,.
\]
But up to isomorphism, the only non-abelian simple group of order $< 96$ is $A_5$, which we know has $m=3$, and this gives a contradiction.

Therefore $m\in\{2,3,4\}$. It is well documented that finite simple groups have no irreducible representations of degree~$2$ (see e.g.~\cite[Corollary 22.13]{JamLie} for an elementary proof), and it follows from Theorem \ref{t:Blichfeldt 4-dim} that $m\neq 4$. 
The only remaining possibility is that $m=3$.
By Theorem \ref{t:Blichfeldt 3-dim}, this implies that $S\cong A_5$ or $S\cong \PSL(2,7)$; and since $\cp(\PSL(2,7)) = \frac{1}{28} < \frac{1}{20}$, the second case is ruled out. We conclude that $S\cong A_5$, as required.
\end{proof}

\end{subsection}

\end{section}

\newpage


\begin{thebibliography}{Cho23b}

\bibitem[Bli05]{Blichfeldt_PGL4}
{\sc H.~F. Blichfeldt}, {\em The finite, discontinuous primitive groups of
  collineations in four variables}, Math. Ann., 60 (1905), pp.~204--231.

\bibitem[Bli07]{Blichfeldt_PGL3}
\leavevmode\vrule height 2pt depth -1.6pt width 23pt, {\em The finite,
  discontinuous, primitive groups of collineations in three variables}, Math.
  Ann., 63 (1907), pp.~552--572.

\bibitem[Cho23a]{YC_minorant}
{\sc Y.~Choi}, {\em An explicit minorant for the amenability constant of the
  {F}ourier algebra}, Int. Math. Res. Not. IMRN,  (2023), pp.~19390--19430.

\bibitem[Cho23b]{MO_monotone}
\leavevmode\vrule height 2pt depth -1.6pt width 23pt, {\em Moments of character
  degrees - is this result new or folklore?}
\newblock MathOverflow, 2023.
\newblock URL:https://mathoverflow.net/q/437938 (version: 2023-01-08).

\bibitem[Gal70]{Gallagher_cp}
{\sc P.~X. Gallagher}, {\em The number of conjugacy classes in a finite group},
  Math. Z., 118 (1970), pp.~175--179.

\bibitem[GR06]{GR_cp}
{\sc R.~M. Guralnick and G.~R. Robinson}, {\em On the commuting probability in
  finite groups}, J. Algebra, 300 (2006), pp.~509--528.

\bibitem[JL01]{JamLie}
{\sc G.~James and M.~Liebeck}, {\em Representations and characters of groups},
  Cambridge University Press, New York, second~ed., 2001.

\bibitem[Joh94]{BEJ_AG}
{\sc B.~E. Johnson}, {\em Non-amenability of the {F}ourier algebra of a compact
  group}, J. London Math. Soc. (2), 50 (1994), pp.~361--374.

\bibitem[KR93]{KovRob}
{\sc L.~G. Kov\'{a}cs and G.~R. Robinson}, {\em On the number of conjugacy
  classes of a finite group}, J. Algebra, 160 (1993), pp.~441--460.

\bibitem[LLW96]{LLW_SM96}
{\sc A.~T.-M. Lau, R.~J. Loy, and G.~A. Willis}, {\em Amenability of {B}anach
  and {$C^*$}-algebras on locally compact groups}, Studia Math., 119 (1996),
  pp.~161--178.

\bibitem[Mil38]{Miller_|G'|=2}
{\sc G.~A. Miller}, {\em Groups whose commutator subgroups are of order
  two}, Amer. J. Math., 60 (1938), pp.~101--106.

\bibitem[San21]{Sanders_CDT-AG}
{\sc T.~Sanders}, {\em Coset decision trees and the {F}ourier algebra}, J.
  Anal. Math., 144 (2021), pp.~227--259.

\bibitem[TV12]{TongViet_large}
{\sc H.~P. Tong-Viet}, {\em On groups with large character degree sums}, Arch.
  Math. (Basel), 99 (2012), pp.~401--405.

\bibitem[Wal72]{Walter_JFA72}
{\sc M.~E. Walter}, {\em {$W^{\ast} $}-algebras and nonabelian harmonic
  analysis}, J. Functional Analysis, 11 (1972), pp.~17--38.

\end{thebibliography}


\vfill

\newcommand{\address}[1]{{\small\sc#1.}}
\newcommand{\email}[1]{\texttt{#1}}

\noindent
\address{Yemon Choi,
School of Mathematical Sciences,
Lancaster University,
Lancaster LA1 4YF, United Kingdom} 

\noindent
\email{y.choi1@lancaster.ac.uk},\;  \email{y.choi.97@cantab.net}

\end{document}